\documentclass[11pt]{amsart}

\title[]{Landau-Ginzburg Models, D-branes, and Mirror Symmetry}
\author[]{Dmitri Orlov}
\date{}

\address{ Algebra Section, Steklov Mathematical Institute RAN,
Gubkin str. 8, Moscow 119991, RUSSIA}

\email{orlov@mi.ras.ru}

\thanks{
This work was partially supported by  RFBR grants 10-01-93113, 11-01-00336, 11-01-00568, NSh grant 4713.2010.1,
by AG Laboratory HSE, RF government grant, ag. 11.G34.31.0023, and by Simons Center for Geometry and Physics.}

\usepackage{cmap,microtype}

\usepackage{amsmath}
\usepackage{latexsym}
\usepackage{amssymb}
\usepackage{amscd}

\usepackage{amsfonts}
\usepackage[all]{xy}

\newcounter{tmp}

\newtheorem{proposition}{Proposition}[section]

\newtheorem{lemma}[proposition]{Lemma}
\newtheorem{definition}[proposition]{Definition}
\newtheorem{theorem}[proposition]{Theorem}

\newtheorem{conjecture}[proposition]{Conjecture}

\newtheorem{corollary}[proposition]{Corollary}

\newtheorem{example}[proposition]{Example}

\newtheorem{lemma-definition}[proposition]{Lemma-Definition}

\numberwithin{equation}{section}
\DeclareMathOperator{\id}{id}

\DeclareMathOperator{\Mod}{Mod}
\DeclareMathOperator{\End}{End}
\DeclareMathOperator{\md}{mod}

\DeclareMathOperator{\Lag}{Lag}

\renewcommand{\Im}{\mathfrak{Im}}

\newcommand{\Qcoh}{\operatorname{Qcoh}}
\newcommand{\coh}{\operatorname{coh}}

\def\th#1{\par\vspace{0.25cm}
\par\noindent\begingroup \it
\leftskip=0em\hspace{0em}{\bf #1  }}

\def\eth{\par\endgroup}

\def\lto{\longrightarrow}

\newcommand{\bD}{\mathbf{D}}
\newcommand{\bH}{\mathbf{H}}
\newcommand{\bL}{\mathbf{L}}
\newcommand{\bR}{\mathbf{R}}

\def\Perf{\mathfrak{Perf}}
\def\Spec{\operatorname{Spec}}

\def\Sing{\operatorname{Sing}}

\def\Ker{\operatorname{Ker}\,}
\def\Coker{\operatorname{Coker}\,}
\def\Im{\operatorname{Im}\,}

\def\wX{\mathfrak{X}}

\def\FS{\mathfrak{FS}}

\def\dg{\mathfrak{DG}}
\def\dt{{\rm H}^0\dg}

\def\abe{\mathfrak{Pair}}
\def\MF{\mathfrak{MF}}

\def\bE{{\mathbb E}}
\def\bF{{\mathbb F}}
\def\bG{{\mathbb G}}
\def\bT{{\mathbb T}}
\def\ZZ{{\mathbb Z}}
\def\AA{{\mathbb A}}

\def\kk{{\mathbf k}}

\def\ul{\underline}
\def\Hom{{\mathrm H}{\mathrm o}{\mathrm m}}
\def\Hhom{{\mathbb H}{\mathrm o}{\mathrm m}}
\def\Cone{{\ul{\mathbb{C}{\mathrm{one}}}}}

\def\rH{{\operatorname{H}}}

\def\rZ{{\operatorname{Z}}}
\def\Ac{\mathfrak{Ac}}

\begin{document}
\begin{abstract}
This paper is an introduction to D-branes in Landau-Ginzburg models and Homological Mirror Symmetry.
The paper is based on a series of lectures which were given on Second Latin Congress
on Symmetries in Geometry and Physics that took place
at the University of Curitiba, Brazil in December 2010.
\end{abstract}

\maketitle



\section{Triangulated categories}

\subsection{Definition of triangulated categories and their localizations}

Triangulated categories appeared in algebra and geometry as generalization and formalization of notions of derived and homotopy categories.
If we are given an abelian category $\mathcal{A},$ we can consider its
homotopy category $\bH^*(\mathcal{A})$ and then pass to its derived category
$\bD^*(\mathcal{A}),$ where $*$ is reserved for $\{b, +,-, \emptyset\}$ depending which complexes we consider:
 bounded, bounded above, bounded below, or unbounded. The objects of $\bH^*(\mathcal{A})$ are appropriate complexes
of objects of $\mathcal{A}$ and morphisms are morphisms of complexes
up to homotopy. The objects of $\bD^*(\mathcal{A})$ are the same as
those of $\bH^*(\mathcal{A}),$ but we have to invert  all
quasi-isomorphisms, i.e. all morphisms of complexes that induce isomorphisms on cohomology.
In other words we can obtain the derived category as a
localization of the homotopy category with respect to the class of all quasi-isomorphisms:
\[
\bD^*(\mathcal{A}) := \bH^*(\mathcal{A})[\text{Quasi}^{-1}]
\]
It is evident that there is a canonical embedding $\mathcal{A} \hookrightarrow \bD^*(\mathcal{A}),$
which sends an object $A\in\mathcal{A}$ to the  complex $\cdots\to 0 \to A \to 0 \to\cdots$ with one nontrivial term in degree 0.
Both categories $\bH^*(\mathcal{A})$ and $\bD^*(\mathcal{A})$ have natural \emph{triangulated structures}.

\begin{definition}[\cite{Verdier}]
Let $\mathcal{T}$ be an additive category.
The structure of {\sf a triangulated category} on $\mathcal{T}$
is defined by giving of the following data:
\begin{list}{\alph{tmp})}%
{\usecounter{tmp}}
\item an additive autoequivalence
$[1]: \mathcal{T}\lto\mathcal{T}$ (it is called a shift functor or a translation functor),
\item a class of exact (or distinguished) triangles:
$$
X\stackrel{u}{\lto}Y\stackrel{v}{\lto}Z\stackrel{w}{\lto}X[1],
$$
\end{list}
which must satisfy the set of axioms Verdier  T1--T4.
\begin{itemize}
\item[T1.]
    \begin{itemize}
    \item[a)] For each object $X$ the triangle
    $X\stackrel{\id}{\lto}X\lto 0\lto X[1]$
 is exact.
    \item[b)] Each triangle isomorphic to an exact triangle
is exact.
    \item[c)] Any morphism $X\stackrel{u}{\lto} Y$ can be included
in an exact triangle\\
$X\stackrel{u}{\lto} Y\stackrel{v}{\lto} Z\stackrel{w}{\lto}
X[1].$
    \end{itemize}
\item[T2.] A triangle
$X\stackrel{u}{\lto} Y\stackrel{v}{\lto} Z\stackrel{w}{\lto} X[1]$
is exact if and only if the triangle\\
$Y\stackrel{v}{\lto} Z\stackrel{w}{\lto}
X[1]\stackrel{-u[1]}{\lto} Y[1]$ is exact.
\item[T3.] For any two exact triangles and two morphisms
$f, g$ the diagram below
$$
\xymatrix{
X\ar[r]^u \ar[d]_f \ar@{}[dr]|
{\Box}&Y\ar[r]^v\ar[d]^g&Z\ar[r]^w\ar@{.>}[d]^h&X[1]
\ar[d]^{f[1]}\\
X'\ar[r]^{u'}&Y'\ar[r]^{v'}&Z'\ar[r]^{w'}&X'[1].
}
$$
can be completed to a morphism of triangles by a morphism $h:Z\to Z'.$
\item[T4.] For each pair of morphisms
$X\stackrel{u}{\lto} Y\stackrel{v}{\lto} Z$ there is a commutative
diagram
$$
\begin{CD}
X  @>{u}>>  Y @>{x}>> Z' @>>> X[1]\\
@| @V{v}VV   @VV{w}V  @| \\
X   @>>>    Z  @>{y}>>     Y'   @>{w'}>>  X[1]\\
@.   @VVV       @VV{t}V  @VV{u[1]}V\\
@. X' @=  X'  @>{r}>>   Y[1]\\
@. @VV{r}V      @VVV   @.      \\
@. Y[1]  @>{x[1]}>>    Z'[1] @.
\end{CD}
$$
where the first two rows and the two central columns are exact triangles.
\end{itemize}
\end{definition}

This definition is useful and has many applications to algebra, geometry, topology, and even physics. In particular, any
homotopy category $\bH^*(\mathcal{A})$ and any derived category $\bD^*(\mathcal{A})$ have natural triangulated structures.

Now we recall the definition of a localization of categories. Let
$\mathcal{C}$ be a category and let $\Sigma$ be a class of morphisms in
$\mathcal{C}.$ It is well-known that there is a large category
$\mathcal{C}[\Sigma^{-1}]$ and a functor $Q:\mathcal{C}\to \mathcal{C}[\Sigma^{-1}]$ which is
universal among the functors making the elements of $\Sigma$
invertible.
The category $\mathcal{C}[\Sigma^{-1}]$ has a good description
if $\Sigma$ is a multiplicative system.

A family of morphisms $\Sigma$ in a category $\mathcal{C}$
is called {\sf a multiplicative system} if it
satisfies the following conditions:
\begin{itemize}
\item[M1.] all identical morphisms $\id_{X}$
belongs to
$\Sigma$;
\item[M2.] the composition of two elements of
$\Sigma$ belong to
$\Sigma$;
\item[M3.] any diagram
$X'\stackrel{s}{\longleftarrow} X\stackrel{u}{\lto} Y,$
with $s\in\Sigma$ can be completed to a commutative square
$$
\xymatrix{
X \ar[r]^{u}\ar[d]_{s} & Y \ar@{-->}[d]^{t}\\
X' \ar@{-->}[r]^{u'}& Y'}
$$
with $t\in \Sigma$ (the same when all arrows reversed);
\item[M4.] for any two morphisms $f,g$ the existence
of $s\in\Sigma$ with
$fs=gs$ is equivalent to the existence of $t\in \Sigma$
with $tf=tg.$
\end{itemize}

If $\Sigma$ is a multiplicative system then  $\mathcal{C}[\Sigma^{-1}]$
has the following description.
The objects of $\mathcal{C}[\Sigma^{-1}]$ are the objects of
$\mathcal{C}.$
The morphisms  from $X$  to $Y$ in
$\mathcal{C}[\Sigma^{-1}]$ are  pairs
$(s,f)$ in $\mathcal{C}$
of the form
 $$
X\stackrel{f}{\lto}Y'\stackrel{s}{\longleftarrow}  Y,
\qquad s\in \Sigma
$$
 modulo the following equivalence relation:
 $(f,s)\sim (g,t)$
iff there is a commutative diagram
$$
\xymatrix{
&Y'\ar[d]&\\
X\ar[ur]^f \ar[r]^{h} \ar[dr]_{g}& Y''' &Y\ar[ul]_s \ar[l]_{r}\ar[dl]^t\\
&Y''\ar[u]}
$$
with $r\in\Sigma.$

The composition of the morphisms
$(f,s)$  and $(g,t)$ is a morphism $(g'f, s't)$
defined from the following diagram, which exists by M3:
$$
\xymatrix{
&&Z''&&\\
&Y'\ar@{-->}[ur]^{g'}&&Z'\ar@{-->}[ul]_{s'}&\\
X\ar[ur]^{f}&&Y\ar[ul]_{s}\ar[ur]^g && Z\ar[ul]_t}.
$$

It can be checked that $\mathcal{C}[\Sigma^{-1}]$ is a category and there
is a quotient functor
$$
Q: \mathcal{C} \lto \mathcal{C}[\Sigma^{-1}], \quad X\mapsto X, f \mapsto (f,1)
$$
which inverts all elements of $\Sigma$
and it is universal in this sense.

Let $\mathcal{D}$ be a triangulated category and $\mathcal{N}\subset \mathcal{D}$
be a full triangulated subcategory.
Denote by $\Sigma(\mathcal{N})$ a class of morphisms $s$ in $\mathcal{D}$
embedding into an exact triangle
$$
X\stackrel{s}{\lto} Y\lto N\lto X[1]
$$
with $N\in \mathcal{N}.$
It can be checked that $\Sigma(N)$ is a multiplicative system.
We define the quotient category
$$
\mathcal{D}/\mathcal{N}:=\mathcal{D}[\Sigma(\mathcal{N})^{-1}].
$$
We endow the category $\mathcal{D}/\mathcal{N}$ with a translation functor
induced by the translation functor in the category $\mathcal{D}.$
\begin{lemma}
The category $\mathcal{D}/\mathcal{N}$ becomes a triangulated category by taking for
exact triangles such that are isomorphic to the images of exact
triangles in $\mathcal{D}.$ The quotient functor $Q:\mathcal{D}\lto \mathcal{D}/\mathcal{N}$
annihilates $\mathcal{N}.$ Moreover, any exact functor $F: \mathcal{D}\lto \mathcal{D}'$ of
triangulated categories for which $F(X)\simeq 0$ when $X\in \mathcal{N}$
factors uniquely through $Q.$
\end{lemma}

\subsection{Main geometric examples}

\begin{example}[Rings and Modules]{\rm
Let $A$ be a ring. We can consider the abelian category $\Mod-A$ of all (right) $A$\!-modules and take the unbounded
derived category $\bD(\Mod-A).$ This is a triangulated category with arbitrary direct sums. We can also consider
a full triangulated subcategory of $\bD(\Mod-A)$ that is called a triangulated category of perfect complexes and consists
of all bounded complexes of projective $A$\!-modules of finite type. We denote it as $\Perf(A).$ It is not a derived category of any abelian category, but it is a derived category of the exact category of projective modules of finite type (see \cite{Keller2} for definition).

If the ring $A$ is noetherian, then we can also consider the bounded derived category $\bD^b(\md-A)$ of (right) $A$\!-modules of finite type.
It contains the triangulated category of perfect complexes $\Perf(A)$ and they are equivalent when $A$ has a finite global dimension.
} \end{example}

\begin{example}[Schemes] {\rm
The most important example of a derived or triangulated category comes
from a given scheme $(X,\mathcal{O}_X).$ In this case it is natural to consider an abelian category
of sheaves of $\mathcal{O}_X$\!-modules and an abelian category of quasi-coherent sheaves $\Qcoh(X).$
If $X$ is noetherian then we can also
consider an abelian category of coherent
sheaves
$\coh(X).$ We usually work with the unbounded derived
category $\bD(\Qcoh X)$ and the bounded derived category $\bD^b(\coh X).$ It is important and convenient  that
$\bD(\Qcoh X)$ has all direct sums.
If $X$ is noetherian then the natural functor from $\bD^b(\coh X)$ to
$\bD(\Qcoh X)$ is fully faithful and realizes  an equivalence of $\bD^b(\coh X)$
with the full subcategory $\bD^{\emptyset,\, b}_{\coh}(\Qcoh X)$
consisting of all cohomologically bounded complexes with coherent cohomology
(\cite{Il}, Ex.II, 2.2.2).

The category of
$\mathcal{O}_X$\!-modules and the category of abelian-group-valued
sheaves on $X$ depend on the topology of $X$ (e.g.\ Zariski,
etale, flat and so on), while the categories $\Qcoh X$ and
$\coh X$ do not depend on topology.

Another very important triangulated category that appears in geometry is so-called a triangulated category
of \emph{perfect
complexes}
$\mathfrak{Perf}X.$ This category  consists of all complexes of $\mathcal{O}_X$\!-modules which are locally quasi-isomorphic to
a bounded complex of locally free $\mathcal{O}_X$\!-modules of finite type.
For noetherian scheme we have inclusions
\[
\mathfrak{Perf}X \subseteq \bD^b(\coh X) \subset \bD(\Qcoh X).
 \]
The category of perfect complexes $\mathfrak{Perf}X$ coincides with the subcategory of compact objects
in $\bD(\Qcoh X)$ \cite{Neeman}, i.e. such objects $C$ for which the functor $\Hom(C, -)$ commutes with arbitrary direct sums.
This means that the category of perfect complexes can be defined in internal terms of the triangulated category
$\bD(\Qcoh X).$ On the other hand one can
also recover $\bD(\Qcoh X)$ starting from the category of perfect complexes $\mathfrak{Perf}X,$ so these
two categories contain essentially the same information and can be considered in general situation for any
quasi-compact and separated scheme. Note also that
$\mathfrak{Perf}X$ is not a derived category, it is merely
triangulated.

The bounded derived category of coherent sheaves $\bD^b(\coh X)$ can not to be defined in general; we
require that $X$ be noetherian, or at least coherent.
However to any noetherian scheme $X$ we can attach three triangulated categories
$\bD(\Qcoh X),$ $\mathfrak{Perf} X,$ and $\bD^b(\coh X).$ Note also that if $X$ is regular then $\bD^b(\coh X)\cong \mathfrak{Perf} X.$
\bigskip

\noindent\textbf{Claim:} {\it Studying $X$ is actually studying these three triangulated categories $\mathfrak{Perf}X,$
 $\bD^b(\coh X),$ and $\bD(\Qcoh X).$}
} \end{example}

\begin{example}[Voevodsky category of motives, \cite{Voevodsky}] {\rm
Voevodsky's construction of the category of geometric motives provides
another very important example of a triangulated category.  Let $\mathfrak{Sm}/k$ denote
the category of smooth schemes of finite type over $k.$ Let
$\mathfrak{SmCor}/k$ be a category with the same objects, but whose
morphisms are \emph{correspondences}, i.e. $\Hom(X,Y)$ is a free abelian
group generated by integral closed subschemes $Z \subseteq X \times Y$ such that
$Z \to X$ is finite and surjective onto a connected component. The category  $\mathfrak{SmCor}/k$ is additive and one has
$[X]\oplus[Y]=[X\coprod Y].$
Its homotopy category $\bH^b(\mathfrak{SmCor}/k)$ has a natural triangulated structure.
Let us consider a minimal thick triangulated subcategory $T\subset \bH^b(\mathfrak{SmCor}/k)$ that contains all objects of the form
\[
T = \begin{cases} \text{a)\quad} [X \times \mathbb{A}^1] \to [X] &\text{($\mathbb{A}^1$\!-homotopy)} \\
   \text{b)\quad} [U\cap V] \to [U]\oplus[V] \to [X], \quad \text{when}\quad X=U\bigcup V& \text{(Mayer-Vietoris)} \end{cases}
   \]
The triangulated category of \emph{effective geometric motives over $k$} is defined as the quotient
\[
\mathfrak{DM}^\text{eff}_\text{gm}(k) := \bH^b(\mathfrak{SmCor}/k) \bigl/ T.
\]

The category  $\mathfrak{DM}^\text{eff}_\text{gm}(k)$ has a tensor structure, and
$[X]\otimes [Y]=[X\times Y].$ In the category $\mathfrak{DM}^\text{eff}_\text{gm}(k)$ there is a distinguished element $\mathbb{L}$ called the
\emph{Lefschetz motive} or \emph{Tate motive} such that
$[\mathbb{P}^1] = [\text{pt}] \oplus \mathbb{L}$ (and more generally,
$[\mathbb{P}^n] = [\text{pt}] \oplus \mathbb{L} \oplus
\mathbb{L}^{\otimes 2} \oplus \dotsb \oplus \mathbb{L}^{\otimes n}$).
The triangulated category $\mathfrak{DM}_\text{gm}(k)$ of \emph{geometric motives over $k$}
is obtained from the effective category by by formally inverting the functor  of tensor product with the Tate motive $\!-\otimes\mathbb{L}.$
\bigskip

\noindent\textbf{Claim:} {\it Studying (co)homology theories of varieties (schemes) over $k$ is actually
studying the triangulated category $\mathfrak{DM}_\text{gm}(k).$}
} \end{example}

\begin{example}[Symplectic geometry] {\rm A third example of triangulated categories is coming from symplectic geometry.
Let $(X,\omega, B)$ be a symplectic manifold with B-field. We can consider a so-called `derived' Fukaya category
$\bD\mathfrak{Fuk}(X,\omega_{\mathbb{C}})$ which is the homotopy category of an
$A_\infty$\!--Fukaya category $\mathfrak{Fuk}(X,\omega)$ (and is not actually a derived category by construction).
The main objects of these categories
are Lagrangian submanifolds $L \subset X$ together with local systems
$U$ on $L,$ while morphisms are the Floer cohomologies between them.
\bigskip

\noindent\textbf{Claim:} {\it Studying a symplectic manifold with a B-field $(X,\omega, B)$  is actually
studying the triangulated Fukaya category $\bD\mathfrak{Fuk}(X,\omega_{\mathcal{C}}).$}
} \end{example}

Derived and triangulated categories appear in many other places of geometry, algebra, and topology.
The most famous is the stable homotopy category that is the homotopy category of the category of symmetric spectra.

Many natural geometric relations can be described as equivalences or fully faithful functors between corresponded triangulated categories.
For example, McKay correspondence, many of birational transformations like blow ups, flips and flops have such properties \cite{Blow, BO1,Bridgeland, BKR}.
It seems that geometric Langlands correspondence can also be formulated as an equivalence between triangulated categories.

\section{Homological mirror symmetry}

Traditionally,  the mirror symmetry phenomenon is studied
in the most relevant for the physics case of Calabi-Yau manifolds. It was observed that
some pairs of such manifolds ("mirror partners"), which naturally appear in mathematical
models for physics string theory, exhibit properties (have numerical invariants) that are
symmetric to each other. There were several attempts to mathematically formalize
this mirror phenomenon, one of which was the celebrated definition of Homological Mirror Symmetry (HMS)
 stated by Kontsevich in 1994 \cite{Kontsevich}.

Kontsevich was also the first to suggest that Homological Mirror Symmetry can be extended
to a much more general setting, by considering Landau-Ginzburg models. At this moment we know examples of mirror symmetry for
different type of varieties and believe that
mirror symmetry in some sense exists for varieties of all types.
\bigskip

\noindent\textbf{Claim:} {\it Mirrors should exist for all types of varieties.}

\subsection{Mirror symmetry for Calabi-Yau manifolds}

In its original formulation, Kontsevich's celebrated Homological Mirror
Symmetry conjecture \cite{Kontsevich} concerns mirror pairs of Calabi-Yau varieties, for which it predicts
an equivalence between the derived category of coherent sheaves of one variety and the
derived Fukaya category of the other.

From physics point of view in this case we consider superconformal quantum field theory that is called a sigma-model.

\begin{definition}
A {\sf geometric input for a sigma-model} is a data of the following form $(X, I, \omega, B)$ that consists of a smooth manifold $X$
with a complex structure $I,$ an (1,1)-K\"ahler form $\omega$ and a
real closed $2$\!-form $B$ which is called a B-field.
For shortness we will also use the following combination $B+i\omega$ and denote it by $\omega_\mathbb{C}.$
\end{definition}

If $X$ is a (weak) Calabi-Yau manifold (i.e. $\Omega^\text{top}_X \cong \mathcal{O}_X$) then  it is believed that to any
data $(X,I,\omega_\mathbb{C})$ we can attach, in a natural way, a so-called
\emph{superconformal quantum field theory (SCQFT) with
N=(2,2) supersymmetry.}

There is  the following heuristic argument supporting this claim.
First of all, to any such data  one can naturally attach
a classical field theory called the $N=(2,2)$ sigma-model. Its Lagrangian is given by an explicit, although rather
complicated, formula. Infinitesimal symmetries of this classical field theory include
two copies of $N=2$
super-Virasoro algebra (with zero central charge). Second,
one can try to quantize this classical field theory while preserving $N=(2,2)$ superconformal invariance
(up to an unavoidable
central extension). The result of the quantization should be an $N=(2,2)$ SCQFT.

Except for a few special cases, it is not known how to quantize
the sigma-model exactly. On the other hand, one has a perturbative
quantization procedure which works when the volume of the
Calabi-Yau is large. That is, if one rescales the metric by a
parameter $t\gg 1,$ $g_{\mu\nu}\to t^2 g_{\mu\nu},$ and considers
the limit $t\to \infty$ (so called large volume limit), then one
can quantize the sigma-model order by order in $1/t$ expansion. It
is believed that the resulting power series in $1/t$ has a
non-zero radius of convergence, and defines an actual $N=(2,2)$ SCQFT.

Super-Virasoro algebra has an interesting involution called the mirror automorphism.
Suppose we have a pair of $N=(2,2)$ SCQFTs and an isomorphism between them that acts as the identity on the
``left-moving'' $N=2$ super-Virasoro, and acts by the mirror
automorphism on the ``right-moving'' $N=2$ super-Virasoro. In this situation we say that these SCQFTs are mirror symmetric.

Mirror Symmetry can also be extended on the case when 2-dimensional world-sheet has boundaries (see
e.g.~\cite{WittenCS}). This generalization leads to the notion of a D-brane, which plays
a very important role in string theory~\cite{Polch}. A D-brane is
a nice boundary condition for the SCQFT.  For example, one can impose Dirichlet boundary
conditions (i.e. vanishing) on some scalar fields which appear in
the Lagrangian.

If the  field theory has some symmetries, it
is reasonable to require the boundary condition to preserve this symmetry. It is not possible to preserve both of them.
In the $N=(2,2)$ case, we have two copies of $N=2$ super-Virasoro, and we may require
the boundary condition to preserve the diagonal $N=2$ super-Virasoro. Such boundary conditions are called D-branes
of type B, or simply B-branes. One can also exploit the existence of the
mirror automorphism  and consider boundary conditions which preserve a different $N=2$ super-Virasoro subalgebra.
The corresponding branes are called D-branes of type A, or simply A-branes.
One can show that the set of A-branes (or B-branes) has the
structure of a category.

To summarize, to any physicist's Calabi-Yau we can attach two
categories: the categories of A-branes and B-branes. One can argue
that the category of A-branes (resp. B-branes) does not depend on
the extended complex (resp. extended symplectic)
moduli~\cite{WittenCS}. It is obvious
that if two Calabi-Yau manifolds are related by a mirror morphism,
then the A-brane category of the first manifold is equivalent to
the B-brane category of the second one, and vice versa. Obviously,
if two $N=(2,2)$ SCQFTs related two Calabi-Yau manifolds are
isomorphic, then the corresponding categories of A-branes (and
B-branes) are simply equivalent.

\th{Claim:}
a)  A category of D-branes of type A for a sigma model $(X,I,\omega_\mathbb{C})$ is
the derived Fukaya category $\bD\mathfrak{Fuk}(X,\omega).$

b) A category of D-branes of type B for a sigma model $(X,I,\omega_\mathbb{C})$ is
the derived category of coherent sheaves $\bD^b(\coh(X,I)).$
\eth

As we said above  mirror symmetry is a some simple relation between SCQFTs that should interchange D-branes of type A and type B.
This fact can be considered as a definition of Homological Mirror Symmetry.

\begin{definition}[Homological Mirror Symmetry]
We say that two sigma-models $(X,I,\omega_\mathbb{C})$ and $(X^\vee,
I^\vee, \omega^\vee_\mathbb{C})$ are {\sf homologically mirror symmetric} if we have
equivalences of triangulated categories
\[
\begin{array}{llll}
a)& \bD^b\bigl(\coh(X, I)\bigr) & \cong & \bD\mathfrak{Fuk}\bigl(X^\vee, \omega^\vee_\mathbb{C}\bigr)\\
b)& \bD\mathfrak{Fuk}\bigl(X, \omega_\mathbb{C}\bigr) & \cong & \bD^b\bigl(\coh(X^\vee, I^\vee)\bigr).
\end{array}
\]
\end{definition}

In some sense we say that two sigma-models are mirror symmetric to each other if the algebraic variety
$(X,I)$ is ``equal'' to the symplectic manifold $(X^\vee,\omega^\vee_\mathbb{C})$ and
the symplectic manifold $(X, \omega_\mathbb{C})$ is ``equal'' to the algebraic variety $(X^{\vee}, I^{\vee}).$
Thus a passing to triangulated categories allows us to compare an algebraic variety and a symplectic manifold.
\bigskip

\noindent\textbf{Claim:} {\it The notion of a triangulated category allows us to
compare objects from different fields of mathematics and physics.}

\subsection{Bounded derived categories of coherent sheaves}

Let $X$ be a smooth (quasi)-projective variety. To any such variety we can attach the bounded derived category of coherent sheaves $\bD^b(\coh X).$
It is natural to ask the following questions:
How much information of $X$ does $\bD^b(\coh X)$ contain? When do two different varieties have equivalent the bounded
derived categories of coherent sheaves?

There is a reconstruction theorem for varieites of general type and Fano varieties.

\begin{theorem}[\cite{BO}]
Let $X$ be a smooth projective variety such that either canonical sheaf $K_X$ or anticanonical sheaf $K_X^{-1}$ is ample. If $X'$ is another algebraic variety
such that $\bD^b(\coh X) \simeq \bD^b(\coh X').$ Then $X' \cong X.$
\end{theorem}
In this case we can reconstruct a variety from the derived category $D^b(\coh X)$ and, moreover, it can be done directly, i.e. there is a procedure
for the reconstruction.

However, it is much more interesting when we can not to reconstruct a variety and have different varieties
 with equivalent bounded derived categories of coherent sheaves. Such varieties are called Fourier-Mukai partners.
 It happens very seldom and every time any such example has lots of geometric senses and meanings.
The first example is due to Mukai and can be considered as a categorical Fourier transform.

\begin{theorem}[Mukai]
Let $A$ be an abelian variety and $\widehat{A}$ be the dual abelian variety then $\bD^b(\coh A) \simeq
\bD^b(\coh\widehat{A})$ and this equivalence is given by Fourier transform, i.e. the functor is isomorphic
to ${\mathbf R}p_{2*}( p^*_1(-)\otimes \mathcal{P}),$ where $\mathcal{P}$ is a Poincare line bundle
on the product $A\times \widehat{A}.$
\end{theorem}

In this case we have an isomorphism between $N=(2,2)$ SCQFTs as well.
More precisely, there is  the following theorem.
\begin{theorem}[\cite{KO1}]
Let $A$ be an abelian variety and $\omega_\mathbb{C}=B+i\omega$ be a flat 2-form.
Suppose that $A'$ is another abelian variety such that $\bD^b(\coh A')\cong\bD^b(\coh A).$
Then there exists a flat 2-form ${\omega'}_\mathbb{C}$ on $A'$ such that $\operatorname{SCQFT}(A',
{\omega'}_\mathbb{C}) \cong \operatorname{SCQFT}(A, {\omega}_\mathbb{C}).$
\end{theorem}

This example shows that from string theory point of view two varieties with equivalent derived categories of coherent sheaves
should produce the same SCQFT's for a suitable choice of 2-forms on them.

There are known  other examples of equivalences between derived categories of coherent sheaves.
The large class of such examples come from birational geometry.
Let $X$ and $X'$ be birational isomorphic varieties. Suppose that for some (and consequently for any) resolution
$$
\xymatrix{
& \widetilde{X}\ar[dl]_(.4){\pi}
\ar[dr]^(.4){\pi'}&\\
X\ar@{-->}[rr]^{fl}&&X'}
$$
of  the birational isomorphism
$X\stackrel{fl}{\dashrightarrow} X'$ we have $\pi^* K_X=\pi^{\prime*} K_{X'}.$ Such a birational transformation is called a generalized flop.

\begin{conjecture}
If $X$ and $X'$ are related to each other by a generilized flop, then $\bD^b(\coh X')\cong\bD^b(\coh X).$
\end{conjecture}

This conjecture proved for simple examples of flops in \cite{BO1} and by T. Bridgeland for any flop in dimension 3 in \cite{Bridgeland}.
It is also proved for a symplectic flop in \cite{Kawamata, Namikawa}.

If two projective varieties $X_1$ and $X_2$ have equivalent the bounded derived categories of coherent sheaves
$\bD^b(\coh X_1)$ and $\bD^b(\coh X_2)$ then we can ask: how to describe such an equivalence? It is proved
in \cite{K3} for smooth projective varieties and in \cite{LO} for any projective varieties $X_1$ and $X_2$ that any equivalence
$F:\bD^b(\coh X_1)\stackrel{\sim}{\to} \bD^b(\coh X_2)$ can be represented by an object on the product. This means that
$F$ is isomorphic to the functor $\Phi_{\mathcal{E}}:=\bR p_{2*} (p^*_{1}(-)\stackrel{\bL}{\otimes} \mathcal{E}),$
where $\mathcal{E}\in \bD^b(\coh(X_1\times X_2))$ and $p_1, p_2$ are the projection of the product on $X_1$ and $X_2$ respectively.
It is important to mention that there is a theorem of Bertrand To\"{e}n which says that any functor between the differential graded categories
of perfect complexes $\Perf_{dg}(X_1)$ and $\Perf_{dg}(X_2)$
are represented by a perfect complex on the product for any quasi-compact and separated schemes $X_1, X_2$ (see \cite{Toen}).

Finally, I would like to formulate another conjecture which states that the bounded derived category of coherent sheaves
keeps almost all information about the variety.

\begin{conjecture}
If $X$ is a quasi-projective variety, then there are only finitely many quasi-projective Fourier-Mukai partners, i.e.
such  $X'$
that $\bD^b(\coh X') \cong \bD^b(\coh X).$
\end{conjecture}

In all cases, like K3 surfaces, abelian varieties when we can describe all Fourier-Mukai partners precisely this conjecture holds.

\subsection{Mirror symmetry for Fano varieties  and  varieties of general type}

If $X$ is not a Calabi-Yau variety then we can not expect that its mirror is a sigma-model. In such cases mirror symmetric object is so called Landau-Ginzburg model.
The most important cases are a Fano variety or a variety of general type (this means that the canonical sheaf $K_X$ is ample or
anti-ample). Let us introduce a notion of a Landau-Ginzburg model that is a generalization of a sigma-model.

\begin{definition}
A {\sf Landau-Ginzburg model} is a collection $(Y, I,\omega, B, W),$
where $Y$ is a smooth variety, $\omega$ is an (1,1)-K\"ahler form, $B$ is a closed real
2-form (B-field), and $W \colon Y \to \mathbb{A}^1$ is a regular function that is called a superpotential.

We also can consider an action of an algebraic group $G$ on $Y$ such that the superpotential $W$ is semi-invariant.
\end{definition}

Note that a sigma-model is a Landau-Ginzburg model with a trivial superpotential $W=0.$
If $(X, I, \omega_{\mathbb{C}})$ is a sigma-model, we have already  defined
\[
DB := \bD^b(\coh X)\quad\text{and}\quad DA := \bD\mathfrak{Fuk}(X,\omega_{\mathcal{C}}).
\]

Now we can formulate Homological Mirror Symmetry relation for LG models.
We say that
two LG models $(Y, I, \omega, B, W)^G$ and $(Y^\vee, I^\vee, \omega^\vee, B^\vee, W^\vee)^{G^{\vee}}$
are mirror-symmetric if there are equivalences
$$
\begin{array}{llll}
a)&  DB(Y, I, W)^G &\cong& DA(Y^\vee, \omega^{\vee}_{\mathbb{C}}, W^\vee)^{G^\vee} \quad \text{ and}\\
b)&  DA(Y, \omega_{\mathbb{C}}, W)^G &\cong& DB(Y^\vee, I^\vee, W^\vee)^{G^\vee}
\end{array}
$$
between categories of $G$\!-equivariant of D-branes of type B and type A in these models.

Thus to talk about HMS for LG models we have to define  the categories of D-branes of type A and  type B in these models.
 The category of D-branes of type B can be defined as categories of matrix factorizations.
They are also closely related to the triangulated categories of singularities.
The categories of D-branes of type A are defined as the categories of vanishing Lagrangian cycles. Constructions of these categories
will be discussed in the next sections.

Now I would like to consider examples of LG models that appear as mirror symmetric models for toric varieties.
In toric case we have a precise procedure of constructing
 of a mirror symmetric LG model. This procedure appeared in papers of
Batyrev, Givental, and Hori--Vafa.

For any toric variety $X$ of dimension $n,$ there is a mirror symmetric LG model $Y$  that is a complex torus $(\mathbb{C}^*)^n$
of the same dimension $n$ with a superpotential $W$  which depends on the fan defining the toric variety $X$ and some parameters
$t_1,\dots, t_k,$ where $k$ is the rank of Picard group of $X.$

We explain the procedure of a construction of the superpotential $W$ considering a few examples.

\begin{example}[Projective space] {\rm The simplest example of a toric variety is the projective space.
Let $X=\mathbb{P}^n$ be the projective space. In the toric
fan we have $e_0+e_1+ \dotsb + e_n=0,$ so we introduce variables $T_i, i=0,\dots,n$ such that $T_0 + T_1 + \cdots + T_n = t,$  where $t$ be a some parameter
that depends on the class of the form $\omega_{\mathbb{C}}$ on the projective space $\mathbb{P}^n.$
Let $Y_i = e^{-T_i}.$ Thus  we have $Y_0 Y_1 \cdots Y_n=e^{-t}.$ Now the superpotential is given by the simple formula
\[
W = \sum_{i=0}^n Y_i= Y_1 + Y_2 + \dotsb + Y_n + \frac{e^{-t}}{Y_1 Y_2 \dotsm Y_n},
\]
that is considered as a function on $(\mathbb{C}^*)^n$ with coordinates $Y_1,\dots, Y_n.$
Thus the mirror symmetric LG model is isomorphic to $(\mathbb{C}^*)^n$ with the superpotential $W$ introduced above
and with a standard exact symplectic form on this complex torus.
}
\end{example}

\begin{example}[Blow up of $\mathbb{P}^{2}$ at 1 point] {\rm Consider the Hirzebruch surface $\mathbb{F}_1$ that is  a blow-up
of $\mathbb{P}^2$ at one point.
We have $T_0 + T_1 + T_2 = t$ and $T_2
+ T' = s,$ and hence $Y_0 Y_1 Y_2 = e^{-t},$ $Y_2 Y'=e^{-s}.$ Now the mirror symmetric LG model is the two-dimensional
complex torus with coordinates $Y_1, Y_2$ and with the following superpotential
\[
W = Y_0+Y_1+Y_2+Y' = Y_1 + Y_2 + \frac{e^{-t}}{Y_1 Y_2} + \frac{e^{-s}}{Y_2}.
 \]
The symplectic form is again the standard exact form on this complex torus.
}
\end{example}

\begin{example}[Blow-ups of $\mathbb{P}^{2}$ at 3 points] {\rm
 Let $S_6$ be the blow-up
of $\mathbb{P}^2$ at three points.
In this case a mirror symmetric LG model is again isomorphic to $(\mathbb{C}^*)^2$ and  the superpotential $W$ has  the following  form
\[
W=Y_1 + Y_2 + Y_1Y_2  + \frac{e^{-r}}{Y_1} + \frac{e^{-s}}{Y_2} + \frac{e^{-t}}{Y_1 Y_2},
\]
where $Y_1, Y_2$ are coordinates on the two-dimensional torus.
}
\end{example}

\section{D-branes of type  B in  Landau-Ginzburg models}

\subsection{Matrix factorisations, affine case}

A mathematical definition of the categories of D-branes of type B in
Landau-Ginzburg models is proposed by M.Kontsevich and it was confirmed
by A.Kapustin and Yu.Li in the paper \cite{KL}.

Suppose we have a Landau-Ginzburg model with a total space $Y$ that is a smooth variety and with
a superpotential $W \colon Y \to \mathbb{A}^1$ that is not constant.
For a
definition of B-branes we don't need a symplectic form on $Y$
which have to be in a LG model too.

For any $\lambda \in \mathbb{A}^1,$
we define a category of \emph{matrix factorizations} denoted by
$\mathfrak{MF}_\lambda(Y,W).$ We
give constructions of this categories under the condition that
$Y=\Spec(A)$ is affine. The general definition that is
more sophisticated see below. Since the category of coherent sheaves on an affine scheme
$Y=\Spec(A)$ is the same as the category of finitely generated
$A$\!-modules we will frequently go from sheaves to modules and
back. Note that under this equivalence locally free sheaves are
the same as projective modules.

Objects of the category $\mathfrak{MF}_\lambda(Y, W)$ are $2$\!-periodic sequences
\[
P_{\bullet}=\dotsb \longrightarrow P_0 \xrightarrow{\ p_0\ } P_1 \xrightarrow{\ p_1\ } P_0 \longrightarrow \dotsb \text{ ,}
\]
where $P_0$ and $P_1$ are projective modules over $A$ and the compositions
${p_0 p_1}$ and ${p_1 p_0}$
are the multiplications with the element $(W-\lambda\cdot \id)\in A,$ i.e.
\[
p_0 p_1=(W-\lambda\cdot \id)\cdot \quad\text{and}\quad p_1 p_0 = (W-\lambda\cdot\id)\cdot
\]
The morphisms between two objects $P_\bullet,$ $Q_\bullet$ are
$2$\!-periodic morphisms up to $2$\!-periodic homotopy.
Thus a morphism $f: P_{\bullet}\to Q_{\bullet}$ in the category $\mathfrak{MF}_\lambda(Y, W)$
 is a pair of morphisms
$f_1: P_1\to Q_1$ and $f_0: P_0\to Q_0$ such that
$f_1 p_0=q_0 f_0$ and $q_1 f_1=f_0 p_1$ modulo $2$\!-homotopy that is by definition
a pair of morphisms
$s: P_0\to Q_1$ and $t:P_1\to Q_0$ such that
$f_1=q_0 t + s p_1$ and $f_0=t p_0 + q_1 s.$

\begin{proposition}[\cite{Trsing}]
The category $\mathfrak{MF}_\lambda(Y, W)$  has a natural structure of a triangulated category for which square
of the shift functor is isomorphic to the identity.
\end{proposition}

It can be shown that the category $\mathfrak{MF}_\lambda(Y, W)$ is trivial
if the fiber over a point $\lambda$ is smooth. All details can be found below when we will discuss
not only affine but a general case.

\begin{definition}
Let $(Y, W)$ be a Landau-Ginzburg model with $Y = \Spec A,$ we define
a {\sf category $DB(Y, W)$ of D-branes of type B  (B-branes)} on $Y$ with the
superpotential $W$ as the product
\[
DB(Y, W) := \prod_{\lambda\in\mathbb{A}^1} \mathfrak{MF}_\lambda(Y, W).
\]
\end{definition}

Since $Y$
is regular, the set of  singular fibers  is finite.
Hence, the category $DB(Y, W)$ of D-branes of type B is a
product of finitely many numbers of triangulated categories.

\subsection{The triangulated category of singularities}

There is another approach to  defining a category of D-branes of type B in Landau-Ginzburg models.
It uses a notion of so called triangulated category of singularities.
Let $X$ be a quasi-projective scheme. Recall that we had defined the
triangulated subcategory of perfect complexes $\mathfrak{Perf}(X) \subseteq
\bD^b(\coh X).$ In the case when $X$ is smooth this inclusion is an equivalence $\mathfrak{Perf}(X)
\cong \bD^b(\coh X).$ However, in the case of a singular scheme the difference between these two categories can be considered as a
measure of singularities of $X.$ It also allows us to define a category of singularities of $X.$

\begin{definition}
The {\sf triangulated category of singularities} of $X$ is defined as the quotient
\[ \bD_\text{sg}(X) := \bD^b(\coh X) \bigl/ \mathfrak{Perf}(X). \]
When we have an action of an algebraic group $G$ on $Y,$ we can also consider an equivariant version of the triangulated category of singularities
\[
\bD_\text{sg}^G(X) := \bD^b(\coh^G X) \bigl/ \mathfrak{Perf}^G(X).
\]
\end{definition}

Triangulated categories of singularities of fibers of the superpotential $W$
in a LG- model have a direct relation to the categories of matrix factorizations defined above.

\begin{theorem}[\cite{Trsing}]
If $(Y, W)$ is a Landau-Ginzburg model and $Y = \Spec A,$ then there is an equivalence
$\mathfrak{MF}_\lambda(Y, W) \cong D_\text{sg}(W^{-1}(\lambda)).$
\end{theorem}

This theorem allows us to recast a definition of the category of D-branes of type B.
\begin{definition} Let $(Y, W)$ be a Landau-Ginzburg model. We can define
a {\sf category $DB(Y, W)$ of D-branes of type B} on $Y$ with the
superpotential $W$ as the product
\[
DB(Y, W) :=
\prod_{\lambda \in \mathbb{A}^1} \bD_\text{sg}(W^{-1}(\lambda)).
\]
\end{definition}
Note that this definition does not require that the total space is affine.

Let us consider the  simplest example. It is a case of the ordinary double point.
\begin{example}{\rm
Let $Y=\mathbb{A}^n$ and $W = \sum_{i=1}^n x_i^2.$ The fiber over zero has the simplest isolated singulary
that is an ordinary double point.
If $n$ is odd, then the triangulated category of singularities  $\bD_{sg}(W^{-1}(0))$ is
equivalent to the category of $\kk$\!-vector spaces with trivial shift functor $[1]=\id.$
It is easy to see that in this case the Grothendieck group $K_0(\bD_{sg}(W^{-1}(0)))$ is isomorphic to
$\mathbb{Z}/2\mathbb{Z}.$
If $n$ is even, then the category $\bD_{sg}(W^{-1}(0))$ is the category of
$\kk$\!-supervector spaces and we have $K_0(\bD_{sg}(W^{-1}(0))) =
\mathbb{Z}.$

The triangulated categories of $A_n$\!-singularities are described in the last section of  \cite{Trsing}.
}
\end{example}

\subsection{Matrix factorizations, general case}

Let us consider a general LG model, the total space of which  is not necessary affine.
By a Landau-Ginzburg model we  mean again the following data: a quasi-projective scheme
$Y$ over a field $\kk$ and a regular function $W$ on $Y$ such that the
morphism $W: Y\to \AA^1_{\kk}$ is flat. (It is equivalent to say
that the map of algebras
$\kk[x]\to \varGamma(\mathcal{O}_Y)$ is an injection.)

Usually in the definition of an LG model we ask that $Y$ be regular.
However for all  considerations regularity of $Y$ is not
necessary. Moreover, it seems that it is very interesting to consider
the case of a  singular total space $Y$ as well.

With any $\kk$\!-point $\lambda\in \AA^1$ we can associate a differential
$\ZZ/2\ZZ$\!-graded category $\dg_{\lambda}(Y, W),$ an exact category
$\abe_{\lambda}(Y, W),$ and a triangulated category $\dt_{\lambda}(Y, W)$ that
is the homotopy category for DG category $\dg_{\lambda}(Y, W).$

Objects of all these categories are 2-periodic sequences of vector bundles
or, in other words,  ordered pairs
$$
\ul{\bE}:=\Bigl(
\xymatrix{
\bE_1 \ar@<0.6ex>[r]^{e_1} &\bE_0 \ar@<0.6ex>[l]^{e_0}
}
\Bigl),
$$ where $\bE_0, \bE_1$ are locally free sheaves of finite type on $Y$
and the compositions ${e_0 e_1}$ and ${e_1 e_0}$ are the
multiplications with the element $(W-\lambda\cdot \id)\in \varGamma (\mathcal{O}_Y).$

Morphisms from $\ul{\bE}$ to $\ul{\bF}$
in the category $\dg_{\lambda}(Y, W)$ form $\ZZ/2\ZZ$\!-graded
complex
$$
{\Hhom}(\ul{\bE},\; \ul{\bF})=\bigoplus_{0\le i,j\le 1} \Hom(\bE_i, \bF_j)
$$
with a natural grading $(i-j)\mod{2},$ and with a
differential $D$ acting on  a homogeneous element $p$
of degree $k$ as
$$
D p=f \cdot p-(-1)^{k}p\cdot e .
$$

The space of morphisms $\Hom(\ul{\bE},\; \ul{\bF})$ in the category
$\abe_{\lambda}(Y, W)$ is the space of morphisms in $\dg_{\lambda}(Y, W)$
which are homogeneous of degree 0 and commute with the differential.

The space of morphisms in the category $\dt_{\lambda}(Y, W)$ is the space of
morphisms in $\abe_{\lambda}(Y, W)$ modulo null-homotopic morphisms, i.e.
$$
\Hom_{\abe_{\lambda}(Y, W)}(\ul{\bE},\; \ul{\bF})= \rZ^0(\Hhom(\ul{\bE},\;
\ul{\bF})),
\qquad
\Hom_{\dt_{\lambda}(Y,W)}(\ul{\bE},\; \ul{\bF})= \rH^0(\Hhom(\ul{\bE},\;
\ul{\bF})).
$$
Thus, a morphism $p:\ul{\bE}\to\ul{\bF}$ in the category
$\abe_{\lambda}(Y,W)$ is a pair of morphisms
$p_1: \bE_1\to \bF_1$ and $p_0: \bE_0\to \bF_0$ such that
$p_1 e_0=f_0 p_0$ and $f_1 p_1=p_0 e_1.$
The morphism $p$ is null-homotopic if there are two morphisms
$s_0: \bE_0\to \bF_1$ and $s_1:\bE_1\to \bF_0$ such that
$p_1=f_0 s_1 + s_0 e_1$ and $p_0=s_1 e_0 + f_1 s_0.$

The category $\dt_{\lambda}(Y, W)$ can be
endowed with a natural structure of a triangulated category.
To specify it we have to define a translation functor $[1]$
and a class of exact triangles.
The translation functor can be defined as a functor
that takes an object $\ul{\bE}$ to the object
$$
\ul{\bE} [1]=\xymatrix{
\Bigl(
\bE_0 \ar@<0.6ex>[r]^{-e_0} &\bE_1 \ar@<0.6ex>[l]^{-e_1}
}
\Bigl),
$$
i.e. it changes the order of the modules and the signs of the maps,
and takes a morphism $p=(p_0, p_1)$ to the morphism
$p[1]=(p_1, p_0).$
We see that the functor $[2]$ is the identity functor.

For any morphism $p:\ul{\bE}\to\ul{\bF}$ from the category
$\abe_{\lambda}(Y, W)$ we define a mapping cone $\Cone(p)$ as an object
$$
\Cone(p)=\xymatrix{ \Bigl( \bF_1\oplus \bE_0 \ar@<0.6ex>[r]^{c_1} &
\bF_0\oplus \bE_1 \ar@<0.6ex>[l]^{c_0} } \Bigl)
$$
such that
$$
c_0=
\begin{pmatrix}
f_0 & p_1\\
0 & -e_1
\end{pmatrix},
\qquad
c_1=
\begin{pmatrix}
f_1 & p_0\\
0 & -e_0
\end{pmatrix}.
$$
There are  maps
$q: \ul{\bF}\to \Cone(p), \; g=(\id , 0)$ and $r: \Cone(p)\to\ul{\bE}[1],\;
r=(0, -\id).$

The standard triangles in
$\dt_{\lambda}(Y, W)$ are defined to be the triangles of the form
$$
\ul{\bE}\stackrel{p}{\lto} \ul{\bF}\stackrel{q}{\lto} \Cone(p)
\stackrel{r}{\lto} \ul{\bE}[1]
$$
for some $p\in \abe_{\lambda}(Y, W).$

\begin{definition}
A triangle
$\ul{\bE} {\to} \ul{\bF} {\to} \ul{\bG}
 {\to} \ul{\bE}[1]$ in $\dt_{\lambda}(Y, W)$
is  called  an {\sf exact triangle} if it is isomorphic to a
standard triangle.
\end{definition}

\begin{proposition}\label{trstr}
The category $\dt_{\lambda}(Y, W)$ endowed with the translation functor
$[1]$ and the above class of exact triangles becomes a triangulated
category.
\end{proposition}

We define a triangulated category $\MF_{\lambda}(Y, W)$ of matrix
factorizations on $(Y, W)$ as a Verdier quotient of $\dt(Y, W)$ by a
triangulated subcategory of ``acyclic'' objects. This quotient will
also be called a triangulated category of D-branes of type B in the
LG model $(Y, W)$ over $\lambda.$
More precisely, for any complex of objects of the category
$\abe_{\lambda}(Y, W)$
$$
\ul{\bE}^{i}\stackrel{d^i}{\lto}
\ul{\bE}^{i+1}\stackrel{d^{i+1}}{\lto}
\cdots\stackrel{d^{j-1}}{\lto}\ul{\bE}^{j}
$$
we can consider a totalization $\ul{\bT}$ of this complex. It is a
pair with
\begin{equation}\label{tot}
\bT_1=\bigoplus_{
k+m\equiv 1\md 2}\bE_k^m,\qquad
\bT_0=\bigoplus_{
k+m\equiv 0 \md 2}\bE_k^m,\quad k=0,1,
\end{equation}
and with $t_l=d^m_k+(-1)^m e_k$ on the component $\bE_k^m,$ where
$l=(k+m)\mod{2}.$

Denote by
$\Ac_{\lambda}(Y, W)$ the minimal full triangulated subcategory that
contains totalizations of all acyclic complexes
in the exact category
$\abe_{\lambda}(Y,W).$
It is easy to see that $\Ac_{\lambda}(Y, W)$ coincides with
the minimal full triangulated subcategory containing  totalizations of
all short exact sequences in
$\abe_{\lambda}(Y,W).$
\begin{definition}[\cite{Nonaff}] We define the {\sf triangulated category of matrix
  factorizations} $\MF_{\lambda}(Y, W)$ on $Y$ with a
superpotential $W$ as the Verdier quotient
$\dt_{\lambda}(Y, W)/\Ac_{\lambda}(Y, W).$
\end{definition}
In particular, this definition implies that any short exact sequence
in $\abe_{\lambda}(Y, W)$
becomes an exact triangle in $\MF_{\lambda}(Y, W).$

With any pair $\ul{\bE}$ on $(Y, W)$ we can associate a short exact sequence
\begin{equation}\label{shseq}
0\lto \bE_1 \stackrel{e_1}{\lto} \bE_0 \lto \Coker e_1 \lto 0
\end{equation}
of coherent sheaves on $Y.$

We can attach to an object $\ul{\bE}$ the sheaf $\Coker e_1.$ This is
a sheaf on $Y.$ But the multiplication with $W$ annihilates it. Hence,
we can consider $\Coker e_1$ as a sheaf on the fiber $W^{-1}(\lambda),$ i.e. there is a
sheaf $\mathcal{E}$ on $W^{-1}(\lambda)$ such that $\Coker e_1\cong i_*\mathcal{E}.$ Any morphism
$p:\ul{\bE}\to\ul{\bF}$ in $\abe_{\lambda}(Y, W)$ gives a morphism between
cokernels.  In this way we get a functor
$\mathrm{Cok}:\abe_{\lambda}(Y,W)\to
\coh(W^{-1}(\lambda)).$

It can be shown that the functor $\mathrm{Cok}: \abe_{\lambda}(Y, W) \to \coh(W^{-1}(\lambda))$
induces  exact functors
$\Pi: \dt_{\lambda}(Y, W)\to \bD_{sg}(W^{-1}(\lambda))$ and $\Sigma:
 \MF_{\lambda}(Y, W)\to \bD_{sg}(W^{-1}(\lambda))$ between triangulated categories.

\begin{theorem}[\cite{Nonaff}]\label{main2} Let $Y$ be a quasi-projective scheme. Then
the natural functor $\Sigma: \MF_{\lambda}(Y, W)\to \bD_{sg}(W^{-1}(\lambda))$ is fully
faithful. Moreover, if $Y$ is regular then the functor $\Sigma$ is an equivalence.
\end{theorem}

\section{Properties of triangulated categories of singularities}

\subsection{Localization and completion}

Let $f: X\to X' $ be a morphism of finite Tor-dimension (for example,
a flat morphism or a regular closed embedding). In this case we have an
inverse image functor $\bL f^*:\bD^b(\coh X')\to \bD^b(\coh X).$ It is
clear that the functor $\bL f^*$ sends
 perfect complexes on $X'$ to  perfect complexes on $X.$
Therefore, the functor $\bL f^*$ induces  an exact functor $\bL
\bar{f}^* :\bD_{sg}(X') \to \bD_{sg}(X).$

A fundamental  property of  triangulated categories of singularities
is a property of locality in Zarisky topology. It says that for any
open embedding $j: U\hookrightarrow X,$ for which $\Sing(X)\subset
U,$ the functor $\bar{j}^*:\bD_{sg}(X)\to \bD_{sg}(U)$ is an
equivalence of triangulated categories \cite{Trsing}.

On the other hand, two analytically isomorphic singularities can
have non-equivalent triangulated categories of singularities. It is easy to see that even
double points given by equations $f=y^2-x^2$ and
$g=y^2-x^2-x^3$ have non-equivalent categories of singularities. The
main reason here is that a triangulated category of singularities
is not necessary idempotent complete. This means
that not for each projector $p: C\to C,\; p^2=p$ there is a decomposition of the form $C=\Ker p\oplus \Im
p.$

For any triangulated category $\mathcal{T}$ we can consider  its so called idempotent completion
 (or Karoubian envelope) $\overline{\mathcal{T}}.$ This
is a category that consists of all kernels of all projectors. It has
a natural structure of a triangulated category and the canonical
functor $\mathcal{T}\to\overline{\mathcal{T}}$ is an exact full embedding.
Moreover, the category $\overline{\mathcal{T}}$ is idempotent complete, i.e.
each idempotent $p: C\to C$ in $\overline{\mathcal{T}}$ arises from a
splitting $\Ker p\oplus\Im p.$ We denote by $\overline{\bD_{sg}{X}}$
the idempotent completion of the triangulated categories of
singularities.

For any closed subscheme $Z\subset X$ we can consider the formal
completion of $X$ along $Z$ as a ringed space $(Z,\,
{\underleftarrow{\lim}} \, \mathcal{O}_{X}/\mathcal{J}^n),$ where $\mathcal{J}$ is
the ideal sheaf corresponding to $Z.$ The formal completion actually
depends only on the closed subset $\operatorname{Supp} Z$ and does not depend on a
scheme structure on $Z.$ We denote by $\wX$ the formal completion of
$X$ along its singularities $\Sing(X).$

\begin{theorem}[\cite{Formcomp}] Let $X$ and $X'$ be two quasi-projective schemes. Assume that the formal completions
$\wX$ and $\wX'$ along singularities are isomorphic. Then the
idempotent completions of the triangulated categories of
singularities $\overline{\bD_{sg}(X)}$ and $\overline{\bD_{sg}(X')}$ are
equivalent.
\end{theorem}

Actually, one can show a little bit more. It is proved
that any object of $\bD_{sg}(X)$ is a direct summand of an object in
its full subcategory $\bD^b_{\Sing(X)}(\coh X)/\Perf_{\Sing(X)}(X),$ where
$\bD^b_{\Sing(X)}(\coh X)$ and $\Perf_{\Sing(X)}(X)$ are subcategories of
$\bD^b(\coh X)$ and $\Perf(X)$ respectively, consisting of complexes
with cohomology supported on $\Sing X.$

\subsection{Reduction of dimension}

There is another type of relations between schemes which give equivalences
for triangulated categories of singularities but under which the quotient  categories
$\bD^b_{\Sing(X)}(\coh X)/\Perf_{\Sing(X)}(X)$ are not necessary equivalent.
It is described in \cite{Singnew}.

Let $S$ be a  noetherian regular scheme. Let $\mathcal{E}$ be a vector bundle
on $S$ of rank $r$ and let $s \in H^0(S, \mathcal{E})$ be a  section. Denote
by $X\subset S$ the zero subscheme of $s.$ Assume that the section
$s$ is regular, i.e. the codimension of the subscheme $X$ in $S$ coincides
with the rank $r.$

Consider the  projective bundles $S'=\mathbb{P}(\mathcal{E}^{\vee})$
and $T=\mathbb{P}(\mathcal{E}^{\vee}|_X),$ where $\mathcal{E}^{\vee}$ is the dual
bundle.
The section $s$ induces a section $s'\in H^0( S', \mathcal{O}_{\mathcal{E}}(1))$ of
the Grothendieck line bundle $\mathcal{O}_{\mathcal{E}}(1)$ on $S'.$ Denote by $Y$ the
divisor on $S'$ defined by the section $s'.$ The natural closed
embedding of $T$ into $S'$ goes through $Y.$
All
schemes defined above can be included in the following commutative
diagram.
$$
\xymatrix{ T
\ar[d]_{p} \ar@{->}[r]^{i}
& Y
 \ar[rd]^{\pi}\ar@{->}[r]^{u} & S'
 \ar[d]^{q}
\\
X \ar@{->}[rr]^{j} &  & S}
$$
Consider the composition functor $\bR i_* p^*: \bD^b(\coh X)\to
\bD^b(\coh Y)$ and denote it by $\Phi_{T}.$
\begin{theorem}[\cite{Singnew}]\label{singnew}
Let schemes $X, Y,$ and $T$ be as above. Then the functor
$$
\Phi_{T}:\bD^b(\coh X)\lto \bD^b(\coh Y)
$$
defined by the formula
$
\Phi_{T}(\cdot)=\bR i_* p^*(\cdot)
$
induces a functor
$$
\phi_{T}: \bD_{sg}(X)\lto \bD_{sg}(Y),
$$
which is an equivalence of triangulated categories.
\end{theorem}
The functor $\Phi_{T}=\bR i_{*} p^*$ has a right adjoint functor
which we denote by $\Phi_{T *}.$ It can be represented as a
composition $\bR p_{*} i^{\flat},$ where $i^{\flat}$ is right
adjoint to $\bR i_{*}.$
It is easy to see that all singularities of $Y$ are concentrated  over the singularities of $X,$ hence the functor
$\Phi_{T*}=\bR p_{*} i^{\flat}$ sends the subcategory $\bD^b_{\Sing(Y)}(\coh Y)$
to the subcategory $\bD^b_{\Sing(X)}(\coh X).$ Therefore, we obtain the following corollary.
\begin{corollary}
The functor $\phi_{T*},$
which realizes an equivalence between the triangulated
categories of singularities of $Y$ and $X,$ gives also a functor
$$
\phi_{T*}^{\prime}: \bD^b_{\Sing(Y)}(\coh Y)/\Perf_{\Sing(Y)}(Y)\lto\bD^b_{\Sing(X)}(\coh
X)/\Perf_{\Sing(X)}(X),
$$
and this functor is fully faithful.
\end{corollary}
Note that the functor
$
\phi_{T*}^{\prime}
$
is not an
equivalence in general.

Theorem \ref{singnew} implies the following application for LG models.
Let $S$ be a smooth quasiprojective variety
and let $f, g\in H^0(S, \mathcal{O}_S)$
be  two regular functions.
Suppose that the zero divisor $D\subset S$ defined by the function $g$
is smooth and the restriction of $f$ on $D$ is not constant.
We can consider $D$ as
a Landau-Ginzburg model with
superpotential
$f_{D}: D\lto \AA^1.$
Another Landau-Ginzburg model
is given
by the smooth variety $T=S\times \AA^1$
and the superpotential $W: T \to \AA^1$
defined by the formula $W=f+xg,$ where $x$ is a coordinate on $\AA^1.$
Denote by $T_{\lambda}$ the fiber of $W$ over the point $\lambda.$
For any $\lambda\in \mathbb{A}^1$ there is an equivalence
$\bD_{sg}\bigl(f_{D}^{-1}(\lambda)\bigr) \cong \bD_{sg}\bigl(W^{-1}(\lambda)\bigr).$
Thus, the categories of D-branes of type B for LG models
$(D, f_D)$ and $(T, W)$ are equivalent.

In particular, we obtain so called Kn\"orrer periodicity that asserts an equivalence of the categories of D-branes of type B
in  LG models
$(Y, f)$ and $\bigl(Y\times \mathbb{A}^2_{\{x,y\} }, W=f+x^2+y^2\bigr).$

\subsection{Sigma-models versus  Landau-Ginzburg models}

We take Landau-Ginzburg models as a generalization of sigma-models. Thus, a sigma model is a particular case of a  Landau-Ginzburg model with
the trivial superpotential $W=0.$ It is known that the category of D-branes of type B for the sigma-model $(X, W=0)$ is the bounded derived category of coherent sheaves $\bD^b(\coh X)$ when $X$ is smooth. On the other hand, we have a definition of D-branes of type B in LG models that uses the notion of the triangulated category of singularities. How to see that these two different definitions give the same answer? Let us study this question.
To do it we have to consider the LG model $(X, W=0)$ with the trivial $\bG_{m}$\!-action.

Let $X$ be any a noetherian scheme and let $\Perf(X)$ be the triangulated category of perfect complexes on $X.$
Consider the trivial action of $\bG_m$ on $X.$ Denote by $\Perf^{\bG_m}(X)$ the triangulated category of equivariant
perfect complexes. In other words, this category is the category of $\ZZ$\!-graded perfect complexes.
Any graded perfect complex $P^{\cdot}$ is a direct  sum of the form $\bigoplus_{k=i}^{j} P^{\cdot}_k.$
Hence, the category $\Perf^{\bG_m}(X)$ has a completely orthogonal decomposition of the form
$$
\Perf^{\bG_m}(X)=\bigoplus_{k\in \ZZ} \Perf(X)_{k},
$$
where $\Perf(X)_0$ is the subcategory with the trivial action of $\bG_m,$ while $\Perf(X)_{k}$ is the category $\Perf(X)_0$ twisted by the
respective character of $\bG_m.$

Consider the constant zero-map to the affine line $\AA^1$ endowed with natural $\bG_m$\!-action.
The fiber of this map over $0$ is $X$ itself. However, since this map is not flat we should take the fiber in derived sense, i.e. as a derived scheme.
Thus, we denote by $\mathcal{X}=(X, \mathcal{P}),$ where $\mathcal{P}$ is a sheaf of DG algebras, the derived cartesian product $\mathbf{0}\mathop{\times}\limits_{\;\; \AA^1}^{\bL} X.$  It is easy to see that $\mathcal{P}$
has only two nontrivial terms
\[
{\mathcal{P}}^0\cong \mathcal{O}_X\in\Perf(X)_0\quad\text{and}\quad {\mathcal{P}}^{-1}\cong \mathcal{O}_{X,-1}\in\Perf(X)_{-1}
\]
with the zero differential. Denote by $\Perf^{\bG_m}(\mathcal{P})$ and $\bD_{\Perf(X)}^{\bG_m}(\mathcal{P})$ triangulated categories of
$\bG_m$\!-equivariant perfect complexes over $\mathcal{P}$ and $\bG_m$\!-equivariant complexes of $\mathcal{P}$\!-modules which are perfect as
$\mathcal{O}_X$\!-modules, respectively.
We have two homomorphisms of sheaves of DG-algebras $e: \mathcal{O}_X\to \mathcal{P}$ and $a:\mathcal{P}\to \mathcal{O}_X$ that induce functors
$e^*: \Perf^{\bG_m}(X)\to \Perf^{\bG_m}(\mathcal{P})$ and $a_*: \Perf^{\bG_m}(X)\to \bD_{\Perf(X)}^{\bG_m}(\mathcal{P}).$ Let us consider restrictions of these functors on the subcategories $\Perf(X)_{k}$ and denote them by $e^*_{k}$ and $a_{*k}$ respectively.

\begin{proposition}
The functors $e^*: \Perf(X)_{k}\to\Perf^{\bG_m}(\mathcal{P})$ are fully faithful for any $k\in\ZZ$ and, moreover, there is a following semi-orthogonal decomposition
$$
\Perf^{\bG_m}(\mathcal{P})=\langle \cdots e^*_{-1}\Perf(X)_{-1}, e^*_{0}\Perf(X)_{0}, e^*_{1}\Perf(X)_{1}, \cdots\rangle
$$
\end{proposition}
Note that this decomposition is not completely orthogonal but it is only semi-orthogonal.

Using the functor $a_{*0}:\Perf_0 \to \Perf^{\bG_m}(\mathcal{P})$ we can also obtain a semi-orthogonal decomposition for
the category $\bD_{\Perf(X)}^{\bG_m}(\mathcal{P}).$
\begin{proposition}
The functors $a_{*k}: \Perf(X)_{k}\to\bD_{\Perf(X)}^{\bG_m}(\mathcal{P})$ are fully faithful for any $k\in\ZZ$
and, moreover, there is a following semi-orthogonal decomposition
$$
\bD_{\Perf(X)}^{\bG_m}(\mathcal{P})=\langle \cdots e^*_{-1}\Perf(X)_{-1}, e^*_{0}\Perf(X)_{0}, a_{*0}\Perf(X)_{0}, e^*_{1}\Perf(X)_{1},
e^*_{2}\Perf(X)_{2}\cdots\rangle
$$
\end{proposition}

These two propositions immediately imply the following theorem.
\begin{theorem}
The quotient category $\bD_{\Perf(X)}^{\bG_m}(\mathcal{P})/\Perf^{\bG_m}(\mathcal{P})$ is equivalent to $\Perf(X).$
\end{theorem}

In the case when $X$ is smooth the quotient category $\bD_{\Perf(X)}^{\bG_m}(\mathcal{P})/\Perf^{\bG_m}(\mathcal{P})$ can be considered as
the category of singularities of the DG scheme $\mathcal{X}=(X, \mathcal{P})$ with the trivial action of the group $\bG_m.$
As we saw above the triangulated category of singularities for DG scheme $\mathcal{X}=(X, \mathcal{P})$
is equivalent to the bounded derived category of coherent sheaves $\bD^b(\coh X)\cong\Perf(X).$

If we consider an LG model $(X, W=0)$ without action of the group $\bG_m$ then it can be shown that the triangulated category
of singularities of the DG scheme $\mathcal{X}$ is equivalent to the derived category of 2-periodic complexes $\bD^{\ZZ/2\ZZ}(\coh X).$

\section{A-side of Landau-Ginzburg models}

\subsection{Category of vanishing cycles}

Let $(Y,I,\omega_\mathbb{C},W)$ be a Landau-Ginzburg model over $\mathbb{C}.$ We define a category of D-branes
$DA(Y, I, \omega_\mathbb{C}, W)$ of type A. It does not depend on the complex structure $I.$

As proposed by Kontsevich \cite{KoEns} and Hori--Iqbal--Vafa \cite{HIV},
the category of A-branes associated with a Landau-Ginzburg model
$W:(Y,\omega_{\mathbb{C}})\to\mathbb{C}$ is a Fukaya-type category which contains not only
compact Lagrangian submanifolds of $Y$ but also certain non-compact
Lagrangians whose ends fiber in a specific way above half-lines in $\mathbb{C}.$
In the case where the critical points of $W$ are isolated and
non-degenerate, this category admits an exceptional collection whose
objects are Lagrangian thimbles associated to the critical points.
Following the formalism introduced by Seidel \cite{Se1,SeBook}, we view
it as the derived category of a finite directed
$A_\infty$\!-category $\Lag_{vc}(Y, \omega_{\mathbb{C}}, W,\{\gamma_i\})$ associated to an ordered
collection of arcs $\{\gamma_i\}.$ The reader is referred to \cite{Se1,SeBook} and to \cite{AKO1, AKO2}
for details.

Consider a symplectic fibration $W:(Y,\omega_{\mathbb{C}})\to \mathbb{C}$ with isolated
non-degenerate critical points, and assume for simplicity that the
critical values $\lambda_0,\dots,\lambda_r$ of $W$ are distinct.
Pick a regular value $\lambda_*$ of $W,$ and choose a collection of
arcs $\gamma_0,\dots,\gamma_r\subset\mathbb{C}$ joining $\lambda_*$ to the
various critical values of $W,$ intersecting each other only at $\lambda_*,$
and ordered in the clockwise direction around $\lambda_*.$ Consider the
horizontal distribution defined by the symplectic form: by parallel
transport along the arc $\gamma_i,$ we obtain a Lagrangian thimble $D_i$
and a vanishing cycle $L_i=\partial D_i\subset \Sigma_*$ (where
$\Sigma_*=W^{-1}(\lambda_*)$). After a small perturbation we can
always assume that the vanishing cycles $L_i$ intersect each other
transversely inside $\Sigma_*.$

\begin{definition}[Seidel]\label{def:fs}
The {\sf directed category of vanishing cycles} $\Lag_{vc}(W,\{\gamma_i\})$ is an
$A_\infty$\!-category over a coefficient ring $R$ with objects
$L_0,\dots,L_r$ corresponding to the vanishing cycles;
the morphisms between
the objects are given by the following rule

$$
\Hom(L_i,L_j)=\begin{cases}
CF^*(L_i,L_j;R)=R^{|L_i\cap L_j|} & \mathrm{if}\ i<j\\
R\cdot id & \mathrm{if}\ i=j\\
0 & \mathrm{if}\ i>j;
\end{cases}$$
and the differential $m_1,$ composition $m_2,$ and higher order
products $m_k$
$$m_k:\Hom(L_{i_0},L_{i_1})\otimes \dots\otimes
\Hom(L_{i_{k-1}},L_{i_k}) \to \Hom(L_{i_0},L_{i_k})[2-k]$$
are defined in terms of Lagrangian Floer homology inside
$\Sigma_*.$
\end{definition}

More precisely, the composition
$m_k$ is trivial when the inequality $i_0<i_1<\dots<i_k$ fails to hold.
When $i_0<\dots<i_k,$ the operations $m_k$ is defined by fixing a generic
$\omega$\!-compatible almost complex structure on $\Sigma_*$ and counting
pseudo-holomorphic maps from a disk with $k+1$ cyclically ordered
marked points on its boundary to $\Sigma_*,$ mapping the marked points
to the given intersection points between vanishing cycles, and the portions
of boundary between them to $L_{i_0},\dots,L_{i_k}$ respectively.

There is a well-defined
$\mathbb{Z}$\!-grading by Maslov index on the Floer complexes $CF^*(L_i,L_j;R)$
once we choose graded Lagrangian lifts of the vanishing cycles.
Considering a nowhere
vanishing  1-form $\Omega\in\Omega^1(\Sigma_*,\mathbb{C})$ and choosing a real
lift of the phase function $\phi_i=\mathrm{arg}(\Omega_{|L_i}):L_i\to S^1$
for each vanishing cycle, one defines a degree of
a given intersection point $p\in L_i\cap L_j$ as
the difference between the phases of $L_i$ and $L_j$ at $p.$

The pseudo-holomorphic disks appearing in
Definition \ref{def:fs} are counted with appropriate weights, and with
signs determined by choices of orientations of the relevant moduli spaces.
The orientation is determined by the choice of a spin structure for each
vanishing cycle $L_i$ (\cite{SeBook}, see also \cite{AKO1}).

The weight attributed to each pseudo-holomorphic map $u$ keeps track
of its relative homology class, which makes it possible to avoid convergence
problems. The usual approach favored by mathematicians is to work over a
Novikov ring, which keeps track of the relative homology class by introducing
suitable formal variables. To remain closer to the physics, we  can use $\mathbb{C}$ as
our coefficient ring, and assign weights according to the symplectic areas.

The weight formula is simplest when there is no B-field; in that case, we
consider untwisted Floer theory, since any flat unitary bundle over the
thimble $D_i$ is trivial and hence restricts to $L_i$ as the trivial bundle.
We then count each map $u:(D^2,\partial D^2)\to (\Sigma_*,\cup L_i)$ with a
coefficient $(-1)^{\nu(u)} \exp(-2\pi\int_{D^2} u^*\omega).$
Hence, given two intersection points $p\in L_i\cap L_j,$ $q\in L_j\cap
L_k$ ($i<j<k$), we have by definition
$$m_2(p,q)=\sum_{\substack{r\in L_i\cap L_k\\\deg r=\deg p+\deg q}}
\Biggl(\sum_{[u]\in \mathcal{M}(p,q,r)}\!\!(-1)^{\nu(u)}
\exp(-2\pi\int_{D^2} u^*\omega) \Biggr)\,r$$
where $\mathcal{M}(p,q,r)$ is the moduli space of pseudo-holomorphic maps
$u$ from the unit disk to $\Sigma_*$ (equipped with a generic $\omega$\!-compatible
almost-complex structure) such that $u(1)=p,$ $u(\mathrm{j})=q,$
$u(\mathrm{j}^2)=r$ (where $\mathrm{j}=\exp(\frac{2i\pi}{3})$), and mapping
the portions of unit circle $[1,\mathrm{j}],$ $[\mathrm{j},\mathrm{j}^2],$
$[\mathrm{j}^2,1]$ to $L_i,$ $L_j$ and $L_k$ respectively.
The other products are defined similarly.

In presence of a B-field, the weights are modified by the fact that we now
consider {\it twisted} Floer homology.
Namely, the weight attributed to a given
pseudo-holomorphic map $u:(D^2,\partial D^2)\to (\Sigma_*,\cup L_i)$
is modified by a factor corresponding to the holonomy along its boundary,
and becomes
$$
(-1)^{\nu(u)}\,\mathrm{hol}(u(\partial D^2))\,\exp(2\pi i\int_{D^2}
u^*(B+i\omega)).
$$
All details can be found  in \cite{SeBook,AKO1,AKO2}.

Although the category $\Lag_{vc}(Y, \omega_{\mathbb{C}}, W,\{\gamma_i \})$ depends on the chosen ordered
collection of arcs $\{\gamma_i\},$ Seidel has obtained the following result
\cite{Se1}.

\begin{theorem}[Seidel]
If the ordered collection $\{\gamma_i\}$ is
replaced by another one $\{\gamma'_i\},$ then the categories
$\Lag_{vc}(Y, \omega_{\mathbb{C}}, W,\{\gamma_i\})$ and $\Lag_{vc}(Y, \omega_{\mathbb{C}}, W, \{\gamma'_i\})$ differ by a sequence of
mutations.
\end{theorem}

Hence, the category naturally associated to the
fibration $W$ is not the finite $A_\infty$\!-category defined above,
but rather an $A_{\infty}$\!-category of twisted complexes over
$\Lag_{vc}(Y, \omega_{\mathbb{C}}, W, \{\gamma_i\})$ that coincides with the $A_{\infty}$\!-category of finite dimensional $A_{\infty}$\!-modules
over $\Lag_{vc}(Y, \omega_{\mathbb{C}}, W, \{\gamma_i\}).$

\begin{definition}[Fukaya, Seidel]
The {\sf Fukaya-Seidel category} $\FS(Y, \omega_{\mathbb{C}}, W)$ is the
$A_\infty$\!-category of twisted complexes over $\Lag_{vc}(X,\omega_{\mathbb{C}}, W, \{\gamma_i\})$ or in this case it is the $A_{\infty}$\!-category of
finite dimensional $A_{\infty}$\!-modules over $\Lag_{vc}(X,\omega_{\mathbb{C}}, W, \{\gamma_i\}).$
\end{definition}

Now the theorem above gives us the following corollary.
\begin{corollary}[Seidel]
The $A_{\infty}$\!-category $\FS(X,\omega,W)$ does not depend on $\{\gamma_i\}.$
\end{corollary}

Now we are ready to define the category of D-branes of type A.
\begin{definition}
Let $(Y, \omega_{\mathbb{C}}, W)$ be a Landau-Ginzburg model. We define
a {\sf category $DA(Y, \omega_{\mathbb{C}}, W)$ of D-branes of type A  (A-branes)} on $Y$ with the
superpotential $W$ as follows
\[
DA(Y, \omega_{\mathbb{C}}, W) :=\bH\FS(Y, \omega_{\mathbb{C}}, W)\cong \bD\FS(Y, \omega_{\mathbb{C}}, W)\cong \bD^b(\md_{fd}-\Lag_{vc}(X,\omega_{\mathbb{C}}, W, \{\gamma_i\})).
\]
\end{definition}
For the first equivalence, recall that the  derived category and the homotopy category
for $A_\infty$\!-categories are the same.

\subsection{Classical generators and mirror symmetry}

When we would like to prove a homological mirror symmetry for some a given pair of models, we have to show that
two categories are equivalent.
How do we check  whether two triangulated categories are equivalent, in general? Roughly speaking, all ways can be divided into two groups.

\begin{description}
\item[Direct way] For $F \colon \mathcal{N} \to \mathcal{M}$ we prove that $F$ is fully faithful, i.e. show that
for any $A, B\in\mathcal{N}$ there is an isomorphism
  $\Hom(A,B) \xrightarrow{\sim} \Hom(FA,FB).$ After that we check that the functor $F$ is essentially surjective on objects.
\item[Indirect way] Use a notion of {\em classical generator} for a triangulated category and we show that two categories have classical genrators with the same DG (or $A_{\infty}$) algebra of endomorphisms.
\end{description}

\begin{definition}
Let $\mathcal{T}$ be a triangulated category. An object $E \in
\mathcal{T}$ is called a {\sf classical generator of $\mathcal{T}$}
if the smallest full triangulated subcategory $\mathcal{U}\subseteq\mathcal{T}$ that
contains $E$ and is closed under taking direct summands (i.e.\ $A \in
\mathcal{U}$ and $A=B\oplus C$ implies $B,C \in \mathcal{U}$)
coincides with $\mathcal{T}.$
\end{definition}

\begin{example}{\rm
Let $A$ be a (right) noetherian algebra. Consider the triangulated category of perfect complexes
$\Perf(A),$ objects of which, by definition, are bounded complexes of (right) projective modules of finite type.
Then $A$ is a classical generator for $\Perf(A).$ Note that the category $\Perf(A)$ is the full subcategory
of the bounded derived category $D^b(\md-A)$ of (right) modules of finite type. They are not necessary equivalent and, hence, $A$ is not a classical generator for $D^b(\md-A).$
} \end{example}
\begin{example}{\rm
If $(E_1, \dots, E_n)$ is a full exceptional collection for a triangulated category $\mathcal{T},$
then the direct sum $\bigoplus_{i=1}^n E_i$ is a classical generator for $\mathcal{T}.$
} \end{example}

If we have an algebra $A$ and a quasi-compact and separated scheme $X$ then the categories $\Perf(A)$ and $\Perf(X)$ can be also described in the internal terms of the unbounded derived categories
$D(\Mod-A)$ and $D(\Qcoh X)$ of all right modules and all quasi-coherent sheaves. They coincides with full subcategories of compact objects \cite{Keller, Neeman}.
Recall that an object $E$ is called compact if the functor $\Hom(E, -)$ commutes with arbitrary
direct sums.
It is proved by Neeman for quasi-compact and separated schemes and by Bondal and Van den Bergh for quasi-compact and quasi-separated schemes
that the triangulated categories of perfect complexes $\Perf(X)$ have  a classical generator.

\begin{theorem}[\cite{Neeman, BV}]
For any quasi-compact, quasi-separated scheme $X$ the triangulated category of perfect complexes  $\mathfrak{Perf}X$ has a classical generator.
\end{theorem}
It is not completely clear how to construct such a generator for a general quasi-compact and quasi-separated scheme. However for a quasi-projective scheme it can be done.

\begin{proposition}[\cite{Gendim}]
If $X$ is a quasi-projective scheme of dimension $n$ and
  $\mathcal{L}$ is a very ample line bundle on $X,$ then the direct sum
  $\bigoplus_{i=0}^n \mathcal{L}^i$ is a classical generator for
  $\mathfrak{Perf}X.$
\end{proposition}

The bounded derived category of coherent sheaves $\bD^b(\coh X)$ for a scheme $X$ of finite type also has a classical generator.

\begin{theorem}[\cite{Rouquier}]
For a scheme $X$ of finite type, $D^b(\coh X)$ has a classical generator that is actually is a strong generator, i.e. it generates the whole category for the finite number of steps.
\end{theorem}

An existing of a classical generator can help us to prove an equivalence between triangulated categories.
Let $\mathcal{T}$ be a triangulated category that has an enhancement, i.e. there is  some pre-triangulated $A_\infty$\!- or
DG-category $\mathcal{U}$ such that $\mathcal{T}$ is an equivalent to the homotopy category $H^0(\mathcal{U}).$ Suppose it has a classical generator $E\in \mathcal{T}=H^0(\mathcal{U}).$
Consider $E$ as an object of the DG-category $\mathcal{U}$
and take a DG-algebra (or $A_{\infty}$\!-algebra)
\[
A := \bR\Hom_\mathcal{U}(E,E).
 \]
Now results of Bernhard Keller give us the following equivalence

\begin{theorem}[\cite{Keller}]\label{Keller} Suppose that $\mathcal{T}=H^0(\mathcal{U})$ is idempotent complete. Then there is an equivalence
$\mathcal{T} \cong \mathfrak{Perf}(A),$ where the category of perfect complexes $\mathfrak{Perf}(A)$ is, by definition, the smallest
triangulated subcategory in the derived category of A-modules $\bD(\Mod-A)$ that contains $A$ and is closed under taking direct summands.
\end{theorem}

Assume now that we have two triangulated categories $\mathcal{T}_1$ and $\mathcal{T}_2$ which are obtained as homotopy categories of
two pre-triangulated DG-categories $\mathcal {U}_1$ and $\mathcal{U}_2$ respectively. Suppose we have two generators
$E_1\in\mathcal{T}_1$ and $E_2\in\mathcal{T}_2$ for which we can check that the DG-algebras of endomorphisms
$A_i=\bR\Hom_{\mathcal{U}_i}(E_i, E_i), \; i=1,2$  are quasi-isomorphic.
Then by theorem \ref{Keller} these triangulated categories are equivalent. Moreover, DG-categories $\mathcal{U}_1$ and $\mathcal{U}_2$ are quasi-equivalent.

Application of this results is a standard way to establish a Homological Mirror Symmetry between categories of D-branes of type A and D-branes of type B in mirror symmetric models: we  choose such generators in both categories and prove that the DG-algebras of endomorphisms are quasi-isomorphic.

\begin{example}{\rm Let us consider the following simple but not trivial example.
Let $X = \mathbb{P}^1$ be the projective line. Consider the sigma-model with $X=\mathbb{P}^1$
as the target space. The category of D-branes of type B is the bounded derived category
of coherent sheaves $\bD^b(\coh \mathbb{P}^1).$
Let us take a mirror symmetric LG model that has a total space $Y = \mathbb{C}^*$ with a superpotential
$W(z) = z+\frac1z.$ We can take a standard symplectic form on $\mathbb{C}^*$ but the category of D-branes of type A
 in this case does not depend on it.
The superpotential $W$ has two critical points. Hence we have two vanishing Lagrangians
$L_1, L_2.$ The smooth fiber consists of two points and intersection of $L_1$ and $L_2$ is exactly these two points.
Thus the category
$\Lag_{vc}(W, \{\gamma \})$ has two objects $L_1, L_2$ and two-dimensional Hom space from the first to the second, i.e.
it is equivalent to
$A:= \bigl\{ \stackrel{1}{\bullet}
   \rightrightarrows \stackrel{2}{\bullet} \bigr\}
$

On the other hand, the bounded derived category of coherent sheaves $\bD^b(\coh \mathbb{P}^1)$ has a full exceptional
collection $\bigl(\mathcal{O}, \mathcal{O}(1)\bigr)$ and hence it is equivalent to the bounded derived category $\bD^b(\md-{\End(\mathcal{O}\oplus\mathcal{O}(1))}).$
Now observe that $A \cong \End\bigl(\mathcal{O}\oplus\mathcal{O}(1)\bigr).$
Thus we obtain
$$
\bD\FS(\mathbb{C}^*, W) \cong \bD^b(\md-A)\cong\bD^b(\coh\mathbb{P}^1).
$$
This gives an equivalence between the category of D-branes of type B for $\mathbb{P}^1$ and
the category of D-branes of type A in the mirror symmetric LG model $(\mathbb{C}^*, W=z+\frac{1}{z})$
}
\end{example}

\subsection{Mirrors for weighted projective planes, del Pezzo sufaces and their noncommutative deformations}

Let $\mathbb{P}^2(a_0, a_1, a_2)$ be
the weighted projective plane  (here $a_0,a_1,a_2$ are coprime
positive integers). It is natural to consider the weighted projective plane as smooth orbifold.
In this case the bounded derived category of coherent sheaves on $\mathbb{P}^2(a_0,a_1,a_2)$
has a full exceptional collection of line bundles
$
\langle \mathcal{O},\dots, \mathcal{O}(a_0+a_1+a_2-1)\rangle.
$
The mirror LG model is the affine hypersurface
$Y=\{y_0^{a_0}y_1^{a_1}y_2^{a_2}=1\}\subset
(\mathbb{C}^*)^3$ equipped with an exact symplectic form $\omega,$ trivial B-field, and
the superpotential $W=y_0+y_1+y_2.$

It is proved in \cite{AKO1} that the bounded derived category of coherent sheaves (B-branes)
on the weighted projective plane $\mathbb{P}^2(a_0,a_1,a_2)$ is
equivalent to the derived category of vanishing Lagrangian cycles (A-branes) on
the affine hypersurface $Y\subset(\mathbb{C}^*)^3$ with an exact symplectic form and the trivial B-field.
Observe that weighted projective planes are rigid in terms of commutative
deformations, but have a one-dimesional toric noncommutative
deformations $\mathbb{P}^2_{\theta}(a_0,a_1,a_2).$

It was showed in the paper \cite{AKO1} that
this mirror correspondence between derived categories can be extended to
the toric noncommutative deformations $\mathbb{P}^2_{\theta}(a_0,a_1,a_2).$ These noncommutative deformations
are related to non-exact variations of the symplectic structure and the B-field
on the mirror LG model $Y.$ Variations of the symplectic structure $\omega_{\mathbb{C}}$ induces a deformation of the derived category of vanishing Lagrangian cycles.
Thus the main theorem says us the following.

\begin{theorem}\cite{AKO1}
Homological Mirror Symmetry holds for  $\mathbb{P}^2(a_0,a_1,a_2)$ and its noncommutative deformations,
i.e. $\bD^b(\coh \mathbb{P}_{\theta}^2(a_0,a_1,a_2))\cong \bD\FS(Y, \omega_{\mathbb{C}}, W).$
\end{theorem}

Given a del Pezzo surface $S_K$ obtained by blowing up $\mathbb{P}^2$ at $k$
points, the mirror Landau-Ginzburg model is an elliptic
fibration $W_k:M_k\to \mathbb{C}$ with $k+3$ nodal singular fibers, which has
the following properties:

\begin{itemize}
\item[(i)] the fibration $W_k$
compactifies to an elliptic fibration $\overline{W}_k$ over $\mathbb{P}^1$ in which the
fiber above infinity consists of $9-k$ rational components;
\item[(ii)] the compactified fibration $\overline{W}_k$ can be obtained
as a deformation of the elliptic fibration $\overline{W}_0:\overline{M}\to\mathbb{P}^1$ which
compactifies the mirror to $\mathbb{P}^2.$
\end{itemize}

Moreover, the manifold $M_k$ is equipped with a
symplectic form $\omega$ and a B-field $B,$ whose cohomology classes are
determined by the set of points $K$ in an explicit manner.

\begin{theorem}\cite{AKO2}
Given a del Pezzo surface $S_K$ obtained by blowing up $\mathbb{P}^2$ at $k$
points, there exists a complexified symplectic form $\omega_{\mathbb{C}}$ on $M_k$
for which $\bD^b(\coh(S_K))\cong \bD\FS(M_k, W_k, \omega_{\mathbb{C}}).$
\end{theorem}

The mirror map, i.e. the
relation between the cohomology class $[\omega_{\mathbb{C}}]=[B+i\omega]\in H^2(M_k,\mathbb{C})$ and
the positions of the blown up points in $\mathbb{P}^2,$ can be described
explicitly (Prop. 5.1  \cite{AKO2}).

On the other hand, not every choice of $[\omega_{\mathbb{C}}]\in H^2(M_k,\mathbb{C})$
yields a category equivalent to the derived category of coherent sheaves
on a del Pezzo surface. There are two reasons for this. First, certain
specific choices of $[\omega_{\mathbb{C}}]$ correspond to deformations of the
complex structure of $X_K$ for which the surface contains a
$\!-2$\!-curve, which causes the anticanonical class to no longer be ample.

More importantly, deformations of the symplectic structure on $M_k$ need
not always correspond to deformations of the complex structure on $S_K$
(observe that $H^2(M_k,\mathbb{C})$ is larger than $H^{1}(S_K,T_{S_K})$). The
additional deformation parameters on the mirror side can however be
interpreted in terms of noncommutative deformations of the del Pezzo
surface $S_K$ (i.e., deformations of the derived category $\bD^b(\coh(X_K))$).
In this context we have the following theorem.
\begin{theorem}\cite{AKO2}
Given any noncommutative deformation of the del Pezzo surface $S_{K,\mu},$
there exists a complexified symplectic form $\omega_{\mathbb{C}}$ on $M_k$
for which the derived category $\bD^b(\coh(S_{K,\mu}))$ is
equivalent to $\bD\FS(M_k, W_k, \omega_{\mathbb{C}}).$ Conversely, for a generic choice of
$[\omega_{\mathbb{C}}]\in H^2(M_k,\mathbb{C}),$ the derived category of Lagrangian
vanishing cycles $\bD\FS(M_k, W_k, \omega_{\mathbb{C}})$ is equivalent to the
derived category of coherent sheaves of a noncommutative deformation of
a del Pezzo surface.
\end{theorem}

The mirror map is again explicit, i.e. the parameters which determine
the noncommutative del Pezzo surface can be read off in a simple manner
from the cohomology class $[B+i\omega].$
The key point in the determination of the mirror map is that
the parameters which determine the composition tensors in
$\bD\FS(W_k)$ can be expressed explicitly in terms of the
cohomology class $[B+i\omega].$ A remarkable feature of these formulas
is that they can be interpreted in terms of  theta functions on a certain
elliptic curve.

\subsection{Homological Mirror Symmetry -- Summary}

Let us summarize in which cases Homological Mirror Symmetry is proved.
It is know for the projective line. It is proved by A.~Polishchuk and E.~Zaslow
for elliptic curves in \cite{PZ}. For K3 quartic surfaces it is proved by Paul Seidel
in \cite{SeK3}. The cases of del Pezzo surfaces, weighted projective planes and their noncommutative deformations are considered
and proved by D.~Auroux, L.~Katzarkov, and D.~Orlov in \cite{AKO1, AKO2}. Mirror Symmetry for toric varieties was discussed and proved
by M.~Abouzaid in paper \cite{AbToric}. Varieties of general type was discussed in paper \cite{KKOY} and Homological Mirror Symmetry for
curves of genus 2 was proved by P.~Seidel in \cite{Segenus2} and for curves of genus greater than 2 by A.~Efimov in \cite{efimov}. M.~Abouzaid
and I.~Smith considered two dimensional complex tori and proved mirror symmetry in case of standard symplectic form in the paper \cite{ASmith}.
The mirror symmetry for abelian varieties is also discussed in the paper \cite{Fukaya}.
In the paper \cite{AAEKO} authors considered punctures spheres and proved HMS.

\end{document}